\documentclass[12pt]{article}
\usepackage{latexsym, amssymb, amsmath, amscd, amsfonts, epsfig, graphicx, colordvi,verbatim,ifpdf}
\usepackage{amsfonts, amsmath, amssymb}
\usepackage{amssymb,amsfonts,amsmath,latexsym,epsfig,cite, psfrag,eepic,color}
\usepackage{amscd,graphics}
\usepackage{latexsym, amssymb,  amsmath,amscd, amsfonts, epsfig, graphicx, colordvi,amsthm}

\usepackage{graphicx}
\usepackage{color}
\usepackage{ifpdf}
\usepackage{fancybox}

\usepackage{float}
\usepackage{cases}

\usepackage{extarrows}

\newtheorem{thm}{Theorem}[section]

\def\qed{\nopagebreak\hfill{\rule{4pt}{7pt}}
\medbreak}

\numberwithin{equation}{section}

\def\qed{\nopagebreak\hfill{\rule{4pt}{7pt}}
\medbreak}

\newlength{\boxedparwidth}
  {\begin{center} \begin{tabular}{|@{\hspace{.315in}}c@{\hspace{.15in}}|}
                  \hline \\ \begin{minipage}[t]{\boxedparwidth}
                  \setlength{\parindent}{.25in}}%
  {\end{minipage} \\ \\ \hline \end{tabular} \end{center}}

\begin{document}
\begin{center}

 { \large\bf Ramanujan-type Congruences for Overpartitions Modulo 16}
\end{center}

\begin{center}
{William Y.C. Chen}$^{a,b}$, {Qing-Hu Hou}$^{a,b}$, {Lisa H. Sun}$^{a,*}$ and {Li Zhang}$^{a}$ \vskip 2mm

   $^{a}$Center for Combinatorics, LPMC-TJKLC\\
   Nankai University, Tianjin 300071, P. R. China\\[6pt]
   $^{b}$Center for Applied Mathematics\\
Tianjin University,  Tianjin 300072, P. R. China\\[6pt]

   \vskip 2mm

    chen@nankai.edu.cn, hou@nankai.edu.cn,\\
    $^{*}$sunhui@nankai.edu.cn, zhangli427@mail.nankai.edu.cn
\end{center}

\vskip 6mm \noindent {\bf Abstract.}
Let $\overline{p}(n)$ denote the number of overpartitions of $n$. Recently, Fortin-Jacob-Mathieu and Hirschhorn-Sellers independently
obtained  2-, 3- and 4-dissections of the generating function for $\overline{p}(n)$ and derived a number of congruences for $\overline{p}(n)$ modulo $4$, $8$ and $64$ including $\overline{p}(5n+2)\equiv 0 \pmod{4}$, $\overline{p}(4n+3)\equiv 0 \pmod{8}$ and $\overline{p}(8n+7)\equiv 0 \pmod{64}$.  By employing dissection techniques, Yao and Xia  obtained   congruences for $\overline{p}(n)$ modulo  $8, 16$ and $32$, such as $\overline{p}(48n+26) \equiv 0 \pmod{8}$, $\overline{p}(24n+17)\equiv 0 \pmod{16}$ and $\overline{p}(72n+69)\equiv 0 \pmod{32}$. In this paper, we give a 16-dissection of the generating function for $\overline{p}(n)$ modulo 16 and we show that   $\overline{p}(16n+14)\equiv0\pmod{16}$ for $n\ge 0$.
Moreover, by using the $2$-adic expansion of the generating function of $\overline{p}(n)$ due to Mahlburg, we obtain that $\overline{p}(\ell^2n+r\ell)\equiv0\pmod{16}$, where $n\ge 0$, $\ell \equiv -1\pmod{8}$ is an odd prime and $r$ is a positive integer with $\ell \nmid r$. In particular, for $\ell=7$, we get $\overline{p}(49n+7)\equiv0\pmod{16}$ and $\overline{p}(49n+14)\equiv0\pmod{16}$ for $n\geq 0$. We also find four congruence  relations: $\overline{p}(4n)\equiv(-1)^n\overline{p}(n) \pmod{16}$ for $n\ge 0$, $\overline{p}(4n)\equiv(-1)^n\overline{p}(n)\pmod{32}$ for $n$ being not a square of an odd positive integer, $\overline{p}(4n)\equiv(-1)^n\overline{p}(n)\pmod{64}$ for $n\not\equiv 1,2,5\pmod{8}$ and $\overline{p}(4n)\equiv(-1)^n\overline{p}(n)\pmod{128}$ for
$n\equiv 0\pmod{4}$.



\section{Introduction}
\allowdisplaybreaks

The objective of this paper is to derive Ramanujan-type congruences for overpartitions modulo $16$, $32$, $64$ and $128$ by constructing a $16$-dissection of the generating function for overpartitions modulo $16$ and applying the $2$-adic expansion due to  Mahlburg \cite{Mahlburg-2004}.

Recall that an overpartition of the nonnegative integer $n$ is a partition of $n$ where the first occurrence of each distinct part may be overlined. We denote the number of overpartitions of $n$ by $\overline{p}(n)$. For example, there are 8 overpartitions of 3:
\[
3,\ \overline{3},\  2 + 1,\  \overline{2} + 1,\  2 + \overline{1},\  \overline{2 }+ \overline{1},\
1 + 1 + 1 ,\  \overline{1} + 1 + 1.
\]
Overpartitions arise in combinatorics \cite{Corteel-Lovejoy-2004}, $q$-series \cite{Corteel-Hitczenko-2004}, symmetric functions \cite{Brenti-1993}, representation theory \cite{Kang-Kwon-2004},
mathematical physics \cite{Fortin-Jacob-Mathieu-2005} and number theory \cite{Lovejoy-2004,Lovejoy-Mallet-2011}, where they are also called standard MacMahon diagrams, joint partitions, jagged partitions or dotted partitions.

As noted by Corteel and Lovejoy \cite{Corteel-Lovejoy-2004},  the generating function of $\overline{p}(n)$ is  given by
\begin{equation}\label{p-phi1}
\sum_{n\geq0}\overline{p}(n)q^n= \frac{(-q;q)_\infty}{(q;q)_\infty},
\end{equation}
where $|q|<1$ and
\[
(a;q)_\infty=\prod_{n=0}^{\infty} (1-aq^n).
\]
The generating function \eqref{p-phi1} can be written in terms of Ramanujan's theta function $\phi(q)$:
\begin{equation}\label{p-phi}
\sum_{n\geq0}\overline{p}(n)q^n=\frac{1}{\phi(-q)},
\end{equation}
where
\begin{align*}
\phi(q)=\sum_{n=-\infty}^{\infty}q^{n^2}.
\end{align*}
  Mahlburg \cite{Mahlburg-2004} showed that the generating function of $\overline{p}(n)$ has the following $2$-adic expansion
\begin{align}\label{2-adic}
\sum_{n\geq0}\overline{p}(n)q^n=1+\sum_{k=1}^{\infty}2^k\sum_{n=1}^{\infty}(-1)^{n+k}c_k(n)q^n,
\end{align}
 where $c_k(n)$ denotes the number of representations of $n$
 as a sum of $k$ squares of positive integers. Employing the above $2$-adic expansion \eqref{2-adic}, Mahlburg \cite{Mahlburg-2004} showed that $\overline{p}(n)\equiv0 \pmod{64}$  for a set of integers of arithmetic density $1$. Moreover, Kim \cite{1-Kim-2008} showed that $\overline{p}(n)\equiv0 \pmod{128}$  for a set of integers of arithmetic density $1$.

 Congruence properties for $\overline{p}(n)$ have been extensively studied, see, for example, \cite{Chen-Xia-2013,Fortin-Jacob-Mathieu-2005,
1-Hirschhorn-Sellers-2005,2-Hirschhorn-Sellers-2005,1-Kim-2008,2-Kim-2009,
Mahlburg-2004,Treneer-2006,Yao-Xia-2013}.
Fortin, Jacob and Mathieu \cite{Fortin-Jacob-Mathieu-2005}, Hirschhorn and Sellers \cite{1-Hirschhorn-Sellers-2005} independently obtained  2-, 3- and 4-dissections of the generating function for $\overline{p}(n)$ and derived a number of Ramanujan-type congruences for $\overline{p}(n)$ modulo $4$, $8$ and $64$, such as
\begin{align}
&\overline{p}(5n+2)\equiv 0 \pmod{4},\nonumber\\[5pt]
&\overline{p}(4n+3)\equiv 0 \pmod{8},\nonumber\\[5pt]
&\overline{p}(8n+7)\equiv 0 \pmod{64}.\label{H-S-2}
\end{align}
Using dissection techniques, Yao and Xia \cite{Yao-Xia-2013} found some congruences for $\overline{p}(n)$ modulo  $8, 16$ and $32$, such as
\begin{align}
&\overline{p}(48n+26) \equiv 0 \pmod{8},\nonumber\\[5pt]
&\overline{p}(24n+17) \equiv 0 \pmod{16},\nonumber\\[5pt]
&\overline{p}(72n+69)\equiv 0 \pmod{32}.\label{Xia--32}
\end{align}

Applying the $2$-adic expansion \eqref{2-adic}
of the generating function for $\overline{p}(n)$, Kim \cite{2-Kim-2009} proved a conjecture of Hirschhorn and Sellers \cite{1-Hirschhorn-Sellers-2005}, that is, if $\ell$ is an odd prime and $r$ is a quadratic nonresidue modulo $\ell$, then
\begin{align*}
\overline{p}(\ell n+r)\equiv\left\{
 \begin{array}{ll}
   0 \pmod{8} &\mbox{if $\ell\equiv\pm1 \pmod{8}$,}\\[5pt]
   0 \pmod{4} &\mbox{if $\ell\equiv\pm3 \pmod{8}$.}
 \end{array}
   \right.
\end{align*}
Moreover, Kim   obtained the following congruence
\begin{align}\label{Kim-8}
\overline{p}(n) \equiv 0 \pmod{8},
\end{align}
where $n$ is neither a square nor twice a square.

It should be noted that Kim's congruences \eqref{Kim-8} contain certain Ramanujan-type congruences for $\overline{p}(n)$ modulo $8$.
Here are some special cases of \eqref{Kim-8}. The detailed proofs
are omitted. For example,
we get the following   Ramanujan-type congruences for $\overline{p}(n)$ modulo $8$.
For $n\geq0$, we have
\begin{align*}
&\overline{p}(8n+5)\equiv 0 \pmod{8}, \\[5pt]
&\overline{p}(8n+6)\equiv 0 \pmod{8},\\[5pt]
&\overline{p}(12n+10)\equiv 0 \pmod{8},\\[5pt]
&\overline{p}(16n+10)\equiv 0 \pmod{8},\\[5pt]
&\overline{p}(20n+6)\equiv 0 \pmod{8}, \\[5pt]
&\overline{p}(20n+14)\equiv 0 \pmod{8}.
\end{align*}

Moreover, as consequences of \eqref{Kim-8},
we obtain three infinite families of Ramanujan-type congruences.
Let $n$ be a nonnegative integer and $\ell$ be an odd prime.
If $r$ is a positive integer with $\ell\nmid r$, then
\begin{align*}
 \overline{p}(\ell^2n+r\ell)\equiv0\pmod{8}.
\end{align*}
If $r$ is an odd positive integer with $\left(\frac{r}{\ell}\right)=-1$, then
\begin{align*}
\overline{p}(2\ell n+r)\equiv0\pmod{8},
\end{align*}
where $\left(\frac{\cdot}{\ell}\right)$ denotes the Legendre symbol.
If $\ell\equiv\pm3 \pmod{8}$ and $\left(\frac{r}{3\ell}\right)=-1$, then
\begin{align*}
\overline{p}(3\ell n+r)\equiv0\pmod{8},
\end{align*}
where $\left(\frac{\cdot}{3\ell}\right)$ denotes the Jocobi symbol.

We are mainly concerned with   congruences for $\overline{p}(n)$ modulo $16$.
We find a $16$-dissection of the generating function for $\overline{p}(n)$ modulo 16, then we establish the following congruence.

\begin{thm}\label{Thm-3}
For $n\geq0$, we have
\begin{align}\label{16n+14--16}
\overline{p}(16n+14)\equiv0\pmod{16}.
\end{align}
\end{thm}

Applying the $2$-adic expansion  \eqref{2-adic}, we derive the following infinite family of congruences for $\overline{p}(n)$ modulo  $16$.

\begin{thm}\label{Thm-6}
Let $n$ be a nonnegative integer, $\ell\equiv-1\pmod{8}$ be an odd  prime and $r$ be a positive integer with $\ell\nmid r$, then we have
\begin{align}\label{l--16}
\overline{p}(\ell^2n+r\ell)\equiv0\pmod{16}.
\end{align}
\end{thm}

For example, when $\ell=7$,  Theorem \ref{Thm-6} gives the following congruences for $n\geq 0$,
\begin{align*}
&\overline{p}(49n+7)\equiv0\pmod{16},\\[5pt]
&\overline{p}(49n+14)\equiv0\pmod{16},\\[5pt]
&\overline{p}(49n+21)\equiv0\pmod{16},\\[5pt]
&\overline{p}(49n+28)\equiv0\pmod{16},\\[5pt]
&\overline{p}(49n+35)\equiv0\pmod{16},\\[5pt]
&\overline{p}(49n+42)\equiv0\pmod{16}.
\end{align*}

The $2$-adic expansion  \eqref{2-adic} can also be used
 to deduce the following congruence relations  $\overline{p}(n)$ modulo $16$, $32$, $64$ and $128$.

\begin{thm}\label{Thm-7}
For  $n\geq0$, we have
\begin{align}\label{4R-mod16}
&\overline{p}(4n)\equiv(-1)^n\overline{p}(n)\pmod{16}.
\end{align}
If $n$ is not  a square of an odd positive integer, then
\begin{align}\label{4R-mod32}
&\overline{p}(4n)\equiv(-1)^n\overline{p}(n)\pmod{32}.
\end{align}
If $n\not\equiv 1,2,5\pmod{8}$, then
\begin{align}\label{4R-mod64}
&\overline{p}(4n)\equiv(-1)^n\overline{p}(n)\pmod{64}.
\end{align}
If $n\equiv 0\pmod{4}$, then
\begin{align}\label{4R-mod128}
&\overline{p}(4n)\equiv(-1)^n\overline{p}(n)\pmod{128}.
\end{align}
\end{thm}

Applying the above congruence relations to the   congruences \eqref{H-S-2}, \eqref{Xia--32} and  \eqref{16n+14--16},  we obtain the following congruences for $n, k\geq 0$,
\begin{align*}
&\overline{p}(4^k(16n+14))\equiv 0 \pmod{16},\\[5pt]
&\overline{p}(4^k(72n+69))\equiv 0 \pmod{32},\\[5pt]
&\overline{p}(4^k(8n+7))\equiv 0 \pmod{64}.
\end{align*}

\section{Proof of Theorem \ref{Thm-3} }

In this section, we obtain a 16-dissection of the generating function for $\overline{p}(n)$ modulo $16$, which implies Theorem \ref{Thm-3}.

Recall that Ramanujan's theta functions $\phi(q)$ and $\psi(q)$ are given by
\begin{align*}
&\phi(q)=\sum_{n=-\infty}^{\infty}q^{n^2},\\[5pt]  &\psi(q)=\sum_{n=0}^{\infty}q^{\frac{n^2+n}{2}}.
\end{align*}

\begin{thm}\label{thm-16dis}
We have a 16-dissection of the generating function for $\overline{p}(n)$ modulo 16, namely,
\begin{align}\label{fz-2}
\sum_{n\geq0}\overline{p}(n)q^n
&\equiv \frac{\phi^{12}}{\phi(-q^{16})^{16}}\Big(\phi^{3}
   - 2\left(4q^{16}\psi^2\psi_2+7\phi^{2}\psi_1\right)q
   + 4\phi\left(q^{16}\psi_2^2+\psi_1^2\right)q^2 \nonumber\\[5pt]
&\hspace{1.5cm}
  + 8\psi_1\left(q^{16}\psi_2^2+\psi_1^2\right)q^3
   + 2\phi\left(4\psi^2+3\phi\psi\right)q^4
   + 8\phi\psi\psi_1q^5 \nonumber\\[5pt]
&\hspace{1.5cm}
   + 8\psi\left(q^{16}\psi_2^2+\psi_1^2\right)q^6
   + 4\phi\psi^2q^8
   - 2\left(7\phi^{2}\psi_2+4\psi^2\psi_1\right)q^9\nonumber\\[5pt]
&\hspace{1.5cm}
   + 8\phi\psi_1\psi_2q^{10}
   + 8\psi_2\left(q^{16}\psi_2^2
   + \psi_1^2\right)q^{11}
   + 8\psi^3q^{12}\nonumber\\[5pt]
&\hspace{1.5cm}
   + 8\phi\psi\psi_2q^{13}\Big) \pmod{16},
\end{align}
where $\phi$, $\psi$, $\psi_1$ and $\psi_2$ denote $\phi(q^{16})$, $\psi(q^{32})$, $\psi_1(q^{16})$ and $\psi_2(q^{16})$, respectively.
\end{thm}

To prove the above formula, we need the following relations
\begin{align}
&\phi(q)=\phi(q^4)+2q\psi(q^8),\label{phi-2-diss}\\[5pt]
&\phi(q)^2=\phi(q^2)^2+4q\psi(q^4)^2,\label{phi^2}\\[5pt]
&\phi(q)\phi(-q)=\phi(-q^2)^2, \label{phi(q)-phi(-q)}
\end{align}
see Berndt \cite[p. 40, Entry 25]{Berndt-1991}.

{\noindent \it Proof of Theorem \ref{thm-16dis}.} We   claim that
\begin{align}\label{phipsi}
\frac{1}{\phi(q)}=
\frac{\phi(-q)\phi(q^2)^2\phi(q^4)^4\phi(q^8)^8}{\phi(-q^{16})^{16}}.
\end{align}
Let $\alpha=e^{\frac{\pi i}{4}}$ and $\beta=e^{\frac{3\pi i}{4}}$.
  Using \eqref{phi(q)-phi(-q)},  we find that
\begin{align*}
\frac{1}{\phi(q)}&=\frac{\phi(-q)\phi(iq)\phi(-iq)\phi(\alpha q)
         \phi(-\alpha q)\phi(\beta q)\phi(-\beta q)}
        {\phi(q)\phi(-q)\phi(iq)\phi(-iq)\phi(\alpha q)
         \phi(-\alpha q)\phi(\beta q)\phi(-\beta q)}\\[5pt]
& = \frac{\phi(-q)\phi(q^2)^2\phi(-iq^2)^2\phi(iq^2)^2}
         {\phi(-q^2)^2\phi(q^2)^2\phi(-iq^2)^2\phi(iq^2)^2} \\[5pt]
&=   \frac{\phi(-q)\phi(q^2)^2\phi(q^4)^4}{\phi(-q^8)^8}\\[5pt]
&= \frac{\phi(-q)\phi(q^2)^2\phi(q^4)^4\phi(q^8)^8}{\phi(-q^8)^8\phi(q^8)^8}\\[5pt]
&=\frac{\phi(-q)\phi(q^2)^2\phi(q^4)^4\phi(q^8)^8}{\phi(-q^{16})^{16}}.
\end{align*}
Therefore, the generating function for $\overline{p}(n)$ can be written as
\[
\sum_{n\geq0}\overline{p}(n)q^n=\frac{1}{\phi(-q)}=\frac{\phi(q)\phi(q^2)^2\phi(q^4)^4
\phi(q^8)^8}{\phi(-q^{16})^{16}}.
\]
Substituting  \eqref{phi-2-diss} and \eqref{phi^2} into the above relation, we obtain that
\begin{align}
&\sum_{n\geq0}\overline{p}(n)q^n\nonumber\\[5pt]
&\quad= \frac{\left(\phi(q^4)+2q\psi(q^8)\right)
    \left(\phi(q^4)^2+4q^2\psi(q^8)^2\right)
    \left(\phi(q^{8})^2+4q^4\psi(q^{16})^2\right)^2
    \left(\phi(q^{16})^2+4q^8\psi(q^{32})^2\right)^4}{\phi(-q^{16})^{16}} \nonumber\\[5pt]
&\quad\equiv \frac{\phi(q^{16})^8}{\phi(-q^{16})^{16}}\Big(\phi(q^4)^3\phi(q^8)^4
+2q\phi(q^4)^2\phi(q^8)^4\psi(q^{8})
+4q^2\phi(q^4)\phi(q^8)^4\psi(q^{8})^2
\nonumber\\[5pt]&\hspace{3.5cm}
+8q^3\phi(q^8)^4\psi(q^{8})^3
+8q^4\phi(q^4)^3\phi(q^8)^2\psi(q^{16})^2
\Big) \pmod{16}.\label{fz}
\end{align}
Write $\psi(q)$ in the following form
\begin{align}\label{psi-2diss}
\psi(q)=\psi_1(q^2)+q\psi_2(q^2),
\end{align}
where
\[
\psi_1(q)=\sum\limits_{n=-\infty}^{\infty}q^{4n^2+n}
\]
and
\[
\psi_2(q)=\sum\limits_{n=-\infty}^{\infty}q^{4n^2-3n}.
\]
Plugging  \eqref{psi-2diss} into \eqref{fz}, we find that
\begin{align}
\sum_{n\geq0}\overline{p}(n)q^n &\equiv \frac{\phi(q^{16})^8}{\phi(-q^{16})^{16}}\Big(\phi(q^4)^3\phi(q^8)^4
+2q\phi(q^4)^2\phi(q^8)^4\big(\psi_1(q^{16})+q^{8}\psi_2(q^{16})\big)
\nonumber\\[5pt]
&\hspace{0.4cm} +4q^2\phi(q^4)\phi(q^8)^4\big(\psi_1(q^{16})+q^{8}\psi_2(q^{16})\big)^2 +8q^3\phi(q^8)^2\big(\psi_1(q^{16})+q^{8}\psi_2(q^{16})\big)^3
\nonumber\\[5pt]
&\hspace{0.4cm}
+8q^4\phi(q^4)^3\phi(q^8)^2\big(\psi_1(q^{32})+q^{16}\psi_2(q^{32})\big)^2
\Big) \pmod{16}.\label{psiz}
\end{align}
Substituting \eqref{phi-2-diss} and \eqref{phi^2} into  \eqref{psiz}, we arrive at \eqref{fz-2}. This
completes the proof. \qed

Notice that  the 16-dissection  of the generating function for $\overline{p}(n)$ modulo 16 given in  Theorem \ref{thm-16dis} contains no terms of  powers of  $q$  congruent to $14$ modulo $16$. So we deduce that $\overline{p}(16n+14)\equiv0\pmod{16}$.

\section{Proof of Theorem \ref{Thm-6}}

In this section, we give a proof of  Theorem \ref{Thm-6} by using the 2-adic expansion \eqref{2-adic} of the generating function for $\overline{p}(n)$.
Recall that Theorem \ref{Thm-6} says that
\begin{align}\label{l-mod 16}
\overline{p}(\ell^2n+r\ell)\equiv0\pmod{16},
\end{align}
where $n\ge 0$, $\ell\equiv-1\pmod{8}$ is an odd  prime and $r$ is a positive integer with $\ell\nmid r$.

\noindent
\emph{Proof of Theorem \ref{Thm-6}}.
By the $2$-adic expansion \eqref{2-adic}, we see that  for $n\geq 1$,
\begin{align*}
\overline{p}(n)\equiv 2(-1)^{n+1}c_1(n)+4(-1)^nc_2(n)+8(-1)^{n+1}c_3(n) \pmod{16}.
\end{align*}
Thus, to prove   congruence \eqref{l-mod 16},
it suffices to show that
\begin{align}
&c_1(\ell^2 n+r\ell) = 0,\label{p^2n+rp--c1--16}\\[5pt]
&c_2(\ell^2 n+r\ell) \equiv 0 \pmod{4},\label{p^2n+rp--c2--16}\\[5pt]
&c_3(\ell^2 n+r\ell) \equiv 0 \pmod{2},\label{p^2n+rp--c3--16}
\end{align}
where $n, \ell$ and $r$  are given in \eqref{l-mod 16}.

First, we consider \eqref{p^2n+rp--c1--16}.  Assume to the contrary that there exists a positive integer $a$ such that
$\ell^2 n+r\ell=a^2$.  It follows that $\ell \mid a^2$.
 Since $\ell$ is a prime, we deduce that $\ell \mid a$, and hence $\ell^2\mid a^2$. Consequently, we have  $\ell \mid r$,
  contradicting the assumption that $\ell\nmid r$. This proves
   \eqref{p^2n+rp--c1--16}.

Next, to prove \eqref{p^2n+rp--c2--16}, it suffices to show that the following equation has no positive integer solution in $x$ and $y$,
\begin{align}\label{l^2n+rl=a^2+b^2--16}
\ell^2 n+r\ell = x^2+y^2.
\end{align}
Otherwise, assume that  $(a,b)$  is a positive integer  solution of \eqref{l^2n+rl=a^2+b^2--16}.
Let $d=\gcd(a,b)$, $a=da_1$ and $b=db_1$.
Then we have
\begin{align*}
\ell^2 n+r\ell = d^2(a_1^2+b_1^2),
\end{align*}
which implies that $\ell\mid d$ or $\ell\mid (a_1^2+b_1^2)$. If $\ell\mid d$, then $\ell^2 \mid d^2$, and hence $\ell\mid r$, which is a contradiction. If $\ell\mid (a_1^2+b_1^2)$, then
\begin{align}\label{l|a_1^2+b_1^2--16}
 a_1^2+b_1^2 \equiv 0 \pmod{\ell}.
\end{align}
Since $\gcd(a_1,b_1)=1$,
 we have $\ell \nmid a_1$ and $\ell \nmid b_1$, so that \eqref{l|a_1^2+b_1^2--16} can be written as
\begin{align*}
\frac{a_1^2}{b_1^2} \equiv -1 \pmod{\ell}.
\end{align*}
This leads to $\left(\frac{-1}{\ell}\right) = 1$, contradicting the assumption $\ell\equiv-1\pmod{8}$. Hence  \eqref{p^2n+rp--c2--16} is proved.

As for \eqref{p^2n+rp--c3--16}, it suffices  to show
that the following equation has an even number of the positive integer solutions in $x,y$ and $z$,
\begin{align}\label{l^2n+rl=x^2+y^2+z^2}
\ell ^2n+r\ell = x^2+y^2+z^2.
\end{align}
Suppose that $(a,b,c)$ is a positive integer solution of \eqref{l^2n+rl=x^2+y^2+z^2}. We consider the following three cases.

Case 1:  $a=b=c$. In this case, \eqref{l^2n+rl=x^2+y^2+z^2} becomes
\begin{align}\label{l^2n+rl=3x^2}
\ell ^2n+r\ell = 3a^2.
\end{align}
Since $\ell \neq 3$ and $\ell\nmid r$,  it is clear that  \eqref{l^2n+rl=3x^2} has no positive integer solution.

Case 2: There are exactly two equal numbers among $a,b$ and $c$. Without loss of generality, we assume that $a=c$, then \eqref{l^2n+rl=x^2+y^2+z^2} becomes
\begin{align}\label{l^2n+rl=2x^2+y^2}
\ell ^2n+r\ell = 2a^2+b^2.
\end{align}
 Using the above argument concerning
 \eqref{l^2n+rl=a^2+b^2--16}, we deduce that  \eqref{l^2n+rl=2x^2+y^2} has no positive integer solution.

Case 3: $a, b$ and $c$ are distinct.
If there exists  a solution $(a,b,c)$, then any permutation of this
triple is also a solution of \eqref{l^2n+rl=x^2+y^2+z^2}. Thus the number of solutions of  \eqref{l^2n+rl=x^2+y^2+z^2} is   even.

In view of the above three cases, we conclude that  \eqref{l^2n+rl=x^2+y^2+z^2} has an even number of positive integer
 solutions, and hence the proof is complete. \qed

\section{Proof of Theorem \ref{Thm-7}}

In this section, we   prove Theorem \ref{Thm-7} by using
the $2$-adic expansion \eqref{2-adic} of the generating function for $\overline{p}(n)$.

\noindent
\emph{Proof of Theorem \ref{Thm-7}}. From the $2$-adic expansion \eqref{2-adic},  we see that for $n\geq 1$ and $k\geq 0$,
\begin{align}\label{2-adic-k}
(-1)^n\overline{p}(n) \equiv -2c_1(n)+2^2c_2(n)+\cdots+(-1)^k2^kc_k(n) \pmod{2^{k+1}}.
\end{align}
Replacing $n$ by $4n$ in \eqref{2-adic-k}, we get
\begin{align}\label{2-adic-4n-k}
\overline{p}(4n) \equiv -2c_1(4n)+2^2c_2(4n)+\cdots+(-1)^k2^kc_k(4n) \pmod{2^{k+1}}.
\end{align}
 By the definition of $c_k(n)$, it is easy to check that for $n\geq 0$,
\[
c_1(n) = c_1(4n), \quad c_2(n) = c_2(4n), \quad c_3(n) = c_3(4n).
\]
Thus \eqref{2-adic-4n-k} becomes
\begin{equation}\label{p4n2}
\overline{p}(4n) \equiv -2c_1(n)+2^2c_2(n)+\cdots+(-1)^k2^kc_k(4n) \pmod{2^{k+1}}.
\end{equation}
Comparing \eqref{2-adic-k} and \eqref{p4n2}, we find that for $n\ge 0$,
\begin{align}\label{2-adic-k-1}
\overline{p}(4n) \equiv &(-1)^n\overline{p}(n)+ 2^4\big(c_4(4n)-c_4(n)\big)\nonumber \\[5pt]
&\qquad +\cdots+(-1)^k2^k\big(c_k(4n)-c_k(n)\big) \pmod{2^{k+1}}.
\end{align}

When $k=3$,  it follows from  \eqref{2-adic-k-1}  that for $n\ge 0$,
\begin{align*}
\overline{p}(4n) \equiv (-1)^n\overline{p}(n) \pmod{16}.
\end{align*}
Setting $k=4$  in \eqref{2-adic-k-1}, we get
\begin{align}
\overline{p}(4n) \equiv (-1)^n\overline{p}(n)+16\big(c_4(4n)-c_4(n)\big) \pmod{32}.\label{n-4n--32}
\end{align}
We claim that
\begin{align}\label{c4--2}
c_4(4n)-c_4(n)\equiv 0 \pmod{2}.
\end{align}
 Observe that the following equation has an even  number of positive integer solutions in $(x,y,z,w)$ such that $x, y, z$ and $w$ are odd,
\begin{align}\label{4n--4sum}
4n=x^2+y^2+z^2+w^2.
\end{align}
 Assume that $(a,b,c,d)$ is a positive integer solution of \eqref{4n--4sum}, where $a, b ,c$ and $d$ are odd.
 Clearly, any permutation of $(a,b,c,d)$ is also a solution of \eqref{4n--4sum}. If there are at least two different numbers among $a$, $b$, $c$ and $d$, then the number of such solutions of  equation \eqref{4n--4sum} is even. Otherwise, we consider the
 case $a=b=c=d$. In this case, we get
$n=a^2$, which contradicts with the assumption that $n$ is not a square of an odd integer. This  proves  \eqref{c4--2}.
Thus, it follows from \eqref{n-4n--32} that
\begin{align*}
\overline{p}(4n) \equiv (-1)^n\overline{p}(n) \pmod{32},
\end{align*}
where $n$ is not a square of an odd positive integer.

To prove \eqref{4R-mod64},
setting  $k=5$ in \eqref{2-adic-k-1}, we find that for $n\geq 0$,
\begin{align}
\overline{p}(4n) \equiv (-1)^n\overline{p}(n)+
16\big(c_4(4n)-c_4(n)\big)-32\big(c_5(4n)-c_5(n)\big) \pmod{64}.\label{n-4n--64}
\end{align}
We claim that for $n\ge 0$,
\begin{align}
&c_4(4n)-c_4(n) \equiv 0 \pmod{4}, \label{4r--c4}\\[5pt]
&c_5(4n)-c_5(n) \equiv 0 \pmod{2}. \label{4r--c5}
\end{align}
As for \eqref{4r--c4}, we need to show that the number of odd positive integer solutions of the equation
\begin{align}\label{4n-4sums-1}
4n = x^2+y^2+z^2+w^2
\end{align}
is a multiple of $4$. Assume that $(a,b,c,d)$ is such a solution of equation \eqref{4n-4sums-1}. If $a=b=c=d$, we get $n=a^2$, which contradicts the assumption $n\not\equiv 1 \pmod{8}$.
If  $a, b, c$ and $d$  are of the pattern $a=b$, $c=d$, but $a\not= c$, regardless of the order, then we get $2n=a^2+c^2$. It contradicts  the assumption that $n\not\equiv 1, 5 \pmod{8}$.  For the other cases,   the number of odd positive solutions of \eqref{4n-4sums-1} is a multiple of $4$. This proves \eqref{4r--c4}.

Congruence \eqref{4r--c5} can be proved by showing that the following equation has an even number of solutions in $(x,y,z,w,v)$,
\begin{align}\label{4n-5sums-1}
4n = x^2+y^2+z^2+w^2+v^2,
\end{align}
where $x, y, z, w, v$ are not all even. Assume that $(a,b,c,d,e)$ is such a solution.
If  $a,b,c,d$ and $e$ are of the pattern $a=c=d=e$, but $a\neq b$,
regardless of the order of $a, b, c, d, e$,  then equation \eqref{4n-5sums-1} becomes
\begin{align}\label{4n--64}
4n = 4a^2+b^2.
\end{align}
Hence $b$ is even. Since $a,b,c,d,e$ are not all even,
  $a$ must be odd. Setting $b=2r$ in \eqref{4n--64}, we deduce that
\begin{align}\label{n--64}
n = a^2+r^2.
\end{align}
But this is impossible, since $n\not\equiv 1, 2, 5 \pmod{8}$. For the other cases, the number of the solutions of equation \eqref{4n-5sums-1} is even. Thus we obtain \eqref{4r--c5}.

Plugging \eqref{4r--c4} and \eqref{4r--c5} into \eqref{n-4n--64}, we deduce that for $n\not\equiv 1,2,5\pmod{8}$,
\begin{align*}
\overline{p}(4n) \equiv (-1)^n\overline{p}(n) \pmod{64}.
\end{align*}
So  \eqref{4R-mod64} is proved.

Setting $k=6$ in \eqref{2-adic-k-1}, we obtain that
\begin{align}\label{n-4n--128}
\overline{p}(4n) \equiv & (-1)^n\overline{p}(n)+
16\big(c_4(4n)-c_4(n)\big)\nonumber \\[5pt]&\quad -32\big(c_5(4n)-c_5(n)\big)+64 \big(c_6(4n)-c_6(n)\big)\pmod{128}.
\end{align}
Using arguments similar to the proofs of congruences \eqref{4r--c4} and \eqref{4r--c5}, we find that for $n\equiv 0 \pmod{4}$,
\begin{align*}
&c_4(4n)-c_4(n) \equiv 0 \pmod{8}, \\[5pt]
&c_5(4n)-c_5(n) \equiv 0 \pmod{4}, \\[5pt]
&c_6(4n)-c_6(n) \equiv 0 \pmod{2}.
\end{align*}
Thus, it follows from \eqref{n-4n--128}  that for $n\equiv 0 \pmod{4}$,
\begin{align*}
\overline{p}(4n) \equiv (-1)^n\overline{p}(n) \pmod{128}.
\end{align*}
This completes the proof.\qed

We remark that congruence  \eqref{4R-mod16} modulo 16 contains the following consequences modulo 4 and 8,
\begin{align*}
&\overline{p}(4n)\equiv \overline{p}(n)\pmod{4},\\[5pt]
&\overline{p}(4n)\equiv (-1)^n\overline{p}(n)\pmod{8},
\end{align*}
which can be used to generate more congruences of $\overline{p}(n)$ modulo $4$ and $8$.

{\noindent \bf Acknowledgments.}  This work was supported by the 973 Project, the PCSIRT Project of the Ministry of Education and the National Science Foundation of China.

\end{document}